\documentclass[10pt]{amsart}
\usepackage{amsfonts}
\usepackage{amsthm}
\usepackage{amsmath}
\usepackage{graphicx}

\title{Van der Corput inequalities for Bessel functions}
\author[\'A. Baricz]{\'Arp\'ad Baricz}
\address{Department of Economics, Babe\c{s}-Bolyai University, 400591 Cluj-Napoca, Romania}
\address{Institute of Applied Mathematics, John von Neumann Faculty of Informatics, \'Obuda University, 1034 Budapest, Hungary}
\email{bariczocsi@yahoo.com}

\author[A. Laforgia]{Andrea Laforgia}
\address{Department of Mathematics, Roma Tre University, Largo San Leonardo Murialdo, 1, 00146 Rome, Italy}
\email{laforgia@mat.uniroma3.it}

\author[T.K. Pog\'any]{Tibor K. Pog\'any}
\address{Faculty of Maritime Studies, University of Rijeka, 51000 Rijeka, Croatia}
\address{Institute of Applied Mathematics, John von Neumann Faculty of Informatics, \'Obuda University, 1034 Budapest, Hungary}
\email{poganj@pfri.hr}

\date{}

\newtheorem{theorem}{Theorem}
\newtheorem{corollary}{Corollary}

\begin{document}
\maketitle

\begin{abstract}
In this note we offer some log-concavity
properties of certain functions related to Bessel functions of the first kind and modified Bessel
functions of the first and second kind, by solving partially a recent conjecture on the log-convexity/log-concavity
properties for modified Bessel functions of the first kind and their derivatives.  Moreover, we give an application of
the mentioned results by extending two inequalities of van der
Corput to Bessel and modified Bessel functions of the first kind. Similar inequalities are proved also for
modified Bessel functions of the second kind, as well as for log-concave probability density functions.
\end{abstract}

\section{\bf Introduction}
\setcounter{equation}{0}

For $\nu$ unrestricted real number the modified Bessel function of
the first kind of order $\nu,$ denoted usually by $I_{\nu},$ is the
particular solution of the differential equation \cite[p.
77]{watson}
\begin{equation}\label{11}
x^2y''(x)+xy'(x)-(x^2+\nu^2)y(x)=0.
\end{equation}
Modified Bessel equation \eqref{11}, which is of frequent occurrence in problems
of mathematical physics, differs from the well-known Bessel's
equation only in the coefficient of $y.$ In fact, the modified
Bessel function $I_{\nu}$ has the infinite series expansion \cite[p.
77]{watson}
\begin{equation}\label{12}
I_{\nu}(x)=\sum_{n\geq
0}\frac{{\left(\frac{x}{2}\right)}^{2n+\nu}}{n!\Gamma(\nu+n+1)}
=\frac{\left(\frac{x}{2}\right)^{\nu}}{\Gamma(\nu+1)}\cdot{}_0F_1\left(\nu+1,\frac{x^2}{4}\right)
\end{equation}
for all $\nu\neq -1,-2,\dots$ and $x\in\mathbb{R},$ where $\Gamma$
denotes the Euler gamma function and
$x\mapsto{}_0F_1\left(\nu+1,\frac{x^2}{4}\right)$ is the hypergeometric function
defined below. We note that the function
$\mathcal{I}_{\nu}:\mathbb{R}\rightarrow[1,\infty),$ defined by
\begin{equation}\label{13}\mathcal{I}_{\nu}(x)={}_0F_1\left(\nu+1,\frac{x^2}{4}\right)=2^{\nu}\Gamma(\nu+1)x^{-\nu}I_{\nu}(x)=
\sum_{n\geq 0}\frac{\left(\frac{1}{4}\right)^n}{{(\nu+1)}_nn!}\cdot
x^{2n},\end{equation} for $\nu\neq -1,-2,\dots,$ (where for $a\neq
0$ the notation ${(a)}_n=a(a+1)(a+2)\dots(a+n-1),$ ${(a)}_0=1$ is
the well-known Pochhammer (or Appell) symbol) is strictly
log-convex, when $\nu>-\frac{1}{2},$ see \cite{neuman}. Observe also that in particular we obtain
$$\mathcal{I}_{-\frac{1}{2}}(x)=\sqrt{\frac{\pi}{2}}\cdot x^{\frac{1}{2}}I_{-\frac{1}{2}}(x)=\cosh
x,$$ and thus actually $\mathcal{I}_{\nu}$ is log-convex for all $\nu\geq
-\frac{1}{2}.$ Moreover, recently it was conjectured that the function $x \mapsto \mathcal{I}_{\nu}(x)$
is strictly log-convex for all $\nu> -1$ and $x \in \mathbb{R},$ see \cite[p. 271]{bane1}. We note that in \cite{GLP},
using H\"older's inequality for integrals, the authors claimed that all
derivatives of the function ${x\mapsto\mathcal{I}_{\nu}(x)}$ are
strictly log-convex for $\nu>-\frac{1}{2}$ and $x\in\mathbb{R}.$ In \cite{MI} it was proved that in fact just the
derivatives of even order of the function $\mathcal{I}_{\nu}$
will be strictly log-convex. Recall that
$\mathcal{I}_{-\frac{1}{2}}(x)=\cosh x,$ thus in this case it is clear
that for all $k\in\{0,1,\dots\}$ the function
$x\mapsto\mathcal{I}^{(2k)}_{-\frac{1}{2}}(x)=\cosh x$ is strictly
log-convex on $\mathbb{R}$ and
$x\mapsto\mathcal{I}^{(2k+1)}_{-\frac{1}{2}}(x)=\sinh x$ is strictly
log-concave on $(0,\infty).$ This suggested the following conjecture (see \cite[p. 101]{bariczPhD}, \cite[p. 187]{barspringer}): the function
$x\mapsto\mathcal{I}^{(2k)}_{\nu}(x)$ is strictly log-convex for all
$k\in\{0,1,\dots\},$ $\nu>-1$ and $x\in\mathbb{R},$ while the function
$x\mapsto\mathcal{I}^{(2k+1)}_{\nu}(x)$ is strictly log-concave for
all $k\in\{0,1,\dots\},$ $\nu>-1$ and $x>0.$

In this paper our aim is to give the full (negative) answer to the above conjectures for the case $k=0.$ We will show that (see the Concluding Remark 2) there exists a $z_{\nu}\geq j_{\nu,1}$ (where $j_{\nu,1}$ is the first positive zero of the Bessel function $J_{\nu}$) such that the function $x \mapsto \mathcal{I}_{\nu}(x)$ is strictly log-convex on $(0,z_{\nu})$ and strictly log-concave on $(z_{\nu},\infty).$ Moreover, we will prove that (see Theorem 1 and its proof) if $\nu\in\left(-1,-\frac{1}{2}\right],$ then the function $x\mapsto\mathcal{I}_{\nu}'(x)$ is strictly log-concave on $(0,\infty),$ while if $\nu>-\frac{1}{2},$ then there exists an $x_{\nu},$ which satisfies $x_{\nu}>\sqrt{2(\nu+3)},$ such that $x\mapsto\mathcal{I}_{\nu}'(x)$ is strictly log-concave on $(0,x_{\nu})$ and strictly log-convex on $(x_{\nu},\infty).$ The paper is organized as follows: in section 2 we study the log-concavity of the function $x\mapsto\mathcal{I}'_{\nu}(x)$ and we apply the main result to extend an inequality of van der Corput concerning hyperbolic functions to modified Bessel functions of the first kind. Moreover, we give the counterparts of the above mentioned results for trigonometric functions and Bessel functions of the first kind. Motivated by these results we present also the counterpart of these results for modified Bessel functions of the second kind. Section 3 contains the proofs of the theorems, and section 4 is devoted for concluding remarks. In section 4, among other things, we point out that the function $x \mapsto \mathcal{I}_{\nu}(x)$
cannot be strictly log-convex on $\mathbb{R}$ for all $\nu> -1,$ and thus we disprove the conjecture from \cite[p. 271]{bane1} (see also \cite[p. 100]{bariczPhD}, \cite[p. 187]{barspringer}).

\section{\bf Inequalities for Bessel and modified Bessel functions}
\setcounter{equation}{0}

\subsection{Modified Bessel functions of the first kind} Our first main result reads as follows.
\begin{theorem}\label{th1}
The following assertions are true:
\begin{enumerate}
\item[(a)] If $\nu\in\left(-1,-\frac{1}{2}\right],$ then the function $x\mapsto\mathcal{I}_{\nu}'(x)$ is strictly log-concave on $(0,\infty).$
\item[(b)] If $\nu>-1,$ then the function $x\mapsto\mathcal{I}_{\nu}'(x)$ is strictly log-concave on $\left(0,\sqrt{2(\nu+3)}\right).$
\item[(c)] If $\nu>-\frac{1}{2},$ then there exists an $x_{\nu},$ which satisfies $x_{\nu}>\sqrt{2(\nu+3)},$ such that $x\mapsto\mathcal{I}_{\nu}'(x)$ is strictly log-concave on $(0,x_{\nu})$ and strictly log-convex on $(x_{\nu},\infty).$
\end{enumerate}
\end{theorem}

The following inequality involving hyperbolic sine and
cosine functions of van der Corput \cite[p. 270]{mi} is of special
interest in this section: \begin{equation}\label{eq.32}|\cosh a
-\cosh b|\geq |a-b|\sqrt{\sinh a \sinh b}\ \ \mbox{for all}\
a,b\geq 0.\end{equation} Taking into account the relation
$\mathcal{I}_{-\frac{1}{2}}(x)=\cosh x$ the inequality (\ref{eq.32}) may
be written as
$$|\mathcal{I}_{-\frac{1}{2}}(a)-\mathcal{I}_{-\frac{1}{2}}(b)|\geq
|a-b|\sqrt{\mathcal{I}_{-\frac{1}{2}}'(a)\mathcal{I}_{-\frac{1}{2}}'(b)}\ \
\mbox{for all}\ a,b\geq 0.$$

The following result is an extension of van der Corput's
inequality (\ref{eq.32}).

\begin{corollary}\label{th2}
If $a,b\geq0$ and $\nu\in\left(-1,-\frac{1}{2}\right],$ then
\begin{equation}\label{eq.34}|\mathcal{I}_{\nu}(a)-\mathcal{I}_{\nu}(b)|\geq
|a-b|\sqrt{\mathcal{I}_{\nu}'(a)\mathcal{I}_{\nu}'(b)}.\end{equation}
Moreover, the above inequality is valid for all $\nu>-1$ and $a,b\in\left(0,\sqrt{2(\nu+3)}\right).$
\end{corollary}

Since for all $n\in\{1,2,\dots\}$ and $\nu\neq -1,-2,\dots$ one has
$(\nu+1){(\nu+2)}_{n-1}={(\nu+1)}_n,$ from (\ref{13}) we have the
following differentiation formula for the function
$\mathcal{I}_{\nu},$ namely
\begin{equation}\label{eq.35}
2(\nu+1)\mathcal{I}_{\nu}'(x)=x\mathcal{I}_{\nu+1}(x),
\end{equation}
where $x\in\mathbb{R}$ and $\nu\neq -1,-2,{\dots}.$ Now taking into
account the relations $\mathcal{I}_{-\frac{1}{2}}(x)=\cosh x$ and
$$\mathcal{I}_{\frac{1}{2}}(x)=\sqrt{\frac{\pi}{2}}\cdot x^{-\frac{1}{2}}I_{\frac{1}{2}}(x)=\frac{\sinh
x}{x},$$ the inequality (\ref{eq.32}) may be written as
$$|\mathcal{I}_{-\frac{1}{2}}(a)-\mathcal{I}_{-\frac{1}{2}}(b)|\geq
|a-b|\sqrt{ab\mathcal{I}_{\frac{1}{2}}(a)\mathcal{I}_{\frac{1}{2}}(b)},$$ and from
(\ref{eq.35}) it follows that the inequality (\ref{eq.34}) is
equivalent to
\begin{equation}\label{eq.36}2(\nu+1)|\mathcal{I}_{\nu}(a)-\mathcal{I}_{\nu}(b)|\geq
|a-b|\sqrt{ab\mathcal{I}_{\nu+1}(a)\mathcal{I}_{\nu+1}(b)}.\end{equation}
We also note that in particular we have
$$\mathcal{I}_{\frac{3}{2}}(x)=3\sqrt{\frac{\pi}{2}}\cdot
x^{-\frac{3}{2}}I_{\frac{3}{2}}(x)=-3\left(\frac{\sinh x}{x^3}-\frac{\cosh
x}{x^2}\right),$$ thus from (\ref{eq.36}) for $\nu=\frac{1}{2}$ we obtain
the following interesting inequality
$$
|b\sinh a-a\sinh b|\geq |a-b|\sqrt{(a\cosh a-\sinh a)(b\cosh
b-\sinh b)} \ \ \mbox{for all}\
a,b\in\left(0,\sqrt{7}\right).
$$

\subsection{Bessel functions of the first kind} Now, let us consider the analogue of $\mathcal{I}_{\nu},$ that is, the function $\mathcal{J}_{\nu}:\mathbb{R}\to (-\infty,1],$ defined by $\mathcal{J}_{\nu}(x)=2^{\nu}\Gamma(\nu+1)x^{-\nu}J_{\nu}(x),$ where $J_{\nu}$ stands for the Bessel function of the first kind. Let us denote by $j_{\nu,n}$ the $n$th positive zero of the Bessel function of the first kind $J_{\nu}.$ The analogue of Theorem \ref{th1} is the following result.

\begin{theorem}\label{th1bessel}
If $\nu>-1$ then the function $x\mapsto-\mathcal{J}_{\nu}'(x)$ is strictly log-concave on $(0,j_{\nu+1,1}).$
\end{theorem}

The following inequality involving sine and
cosine functions of van der Corput \cite[p. 237]{mi} is the analogue of \eqref{eq.32}: \begin{equation}\label{eq.32sine}|\cos a
-\cos b|\geq |a-b|\sqrt{\sin a \sin b}\ \ \mbox{for all}\
a,b\in [0,\pi].\end{equation} Taking into account the relation
$\mathcal{J}_{-\frac{1}{2}}(x)=\cos x$ the inequality (\ref{eq.32sine}) may
be written as
$$|\mathcal{J}_{-\frac{1}{2}}(a)-\mathcal{J}_{-\frac{1}{2}}(b)|\geq
|a-b|\sqrt{\mathcal{J}_{-\frac{1}{2}}'(a)\mathcal{J}_{-\frac{1}{2}}'(b)}\ \
\mbox{for all}\ a,b\in [0,\pi].$$

The following result is an extension of van der Corput's
inequality (\ref{eq.32sine}).

\begin{corollary}\label{th2bessel}
If $a,b\in[0,j_{\nu+1,1}]$ and $\nu>-1,$ then
\begin{equation}\label{eq.34bessel}|\mathcal{J}_{\nu}(a)-\mathcal{J}_{\nu}(b)|\geq
|a-b|\sqrt{\mathcal{J}_{\nu}'(a)\mathcal{J}_{\nu}'(b)}.\end{equation}
\end{corollary}

Observe that $\mathcal{J}_{\nu}$ satisfies the
following differentiation formula
\begin{equation}\label{eq.35bessel}
2(\nu+1)\mathcal{J}_{\nu}'(x)=-x\mathcal{J}_{\nu+1}(x),
\end{equation}
where $x\in\mathbb{R}$ and $\nu\neq -1,-2,{\dots}.$ Now taking into
account the relations
$$\mathcal{J}_{-\frac{1}{2}}(x)=\sqrt{\frac{\pi}{2}}\cdot x^{\frac{1}{2}}J_{-\frac{1}{2}}(x)=\cos
x\ \ \ \mbox{and}\ \ \ \mathcal{J}_{\frac{1}{2}}(x)=\sqrt{\frac{\pi}{2}}\cdot x^{-\frac{1}{2}}J_{\frac{1}{2}}(x)=\frac{\sin
x}{x},$$ the inequality (\ref{eq.32sine}) may be written as
$$|\mathcal{J}_{-\frac{1}{2}}(a)-\mathcal{J}_{-\frac{1}{2}}(b)|\geq
|a-b|\sqrt{ab\mathcal{J}_{\frac{1}{2}}(a)\mathcal{J}_{\frac{1}{2}}(b)},$$ and from
(\ref{eq.35bessel}) it follows that the inequality (\ref{eq.34bessel}) is
equivalent to
\begin{equation}\label{eq.36bessel}2(\nu+1)|\mathcal{J}_{\nu}(a)-\mathcal{J}_{\nu}(b)|\geq
|a-b|\sqrt{ab\mathcal{J}_{\nu+1}(a)\mathcal{J}_{\nu+1}(b)}.\end{equation}
We also note that in particular we have
$$\mathcal{J}_{\frac{3}{2}}(x)=3\sqrt{\frac{\pi}{2}}\cdot
x^{-\frac{3}{2}}J_{\frac{3}{2}}(x)=3\left(\frac{\sin x}{x^3}-\frac{\cos
x}{x^2}\right),$$ thus from (\ref{eq.36bessel}) for $\nu=\frac{1}{2}$ we obtain
the following inequality
$$|b\sin a-a\sin b|\geq |a-b|\sqrt{(a\cos a-\sin a)(b\cos
b-\sin b)} \ \ \mbox{for all}\
a,b\in[0,j_{\frac{3}{2},1}],$$
where $j_{\frac{3}{2},1}=4.493409457\dots$ is the first positive root of the equation $\sin x=x\cos x.$

\subsection{Modified Bessel functions of the second kind} Finally, for $\nu>0$ let us consider the function $\mathcal{K}_{\nu}:(0,\infty)\to(0,\infty),$ defined by $\mathcal{K}_{\nu}(x)=2^{\nu-1}\Gamma(\nu)x^{\nu}K_{\nu}(x),$ where $K_{\nu}$ is the modified Bessel function of the second kind, called sometimes as the MacDonald function. The analogue of Theorem \ref{th1} is the following result.

\begin{theorem}\label{th3k}
If $\nu\geq \frac{3}{2}$ then the function $x\mapsto-\mathcal{K}_{\nu}'(x)$ is strictly log-concave on $(0,\infty).$
\end{theorem}

The following result is an analogue of van der Corput's
inequality (\ref{eq.32}) for modified Bessel functions of the second kind.

\begin{corollary}\label{th4k}
If $a,b\geq0$ and $\nu\geq\frac{3}{2},$ then
\begin{equation}\label{eq.34k}|\mathcal{K}_{\nu}(a)-\mathcal{K}_{\nu}(b)|\geq
|a-b|\sqrt{\mathcal{K}_{\nu}'(a)\mathcal{K}_{\nu}'(b)}.\end{equation}
\end{corollary}

We note that in particular we have
\begin{equation}\label{k3/2}\mathcal{K}_{\frac{3}{2}}(x)=\sqrt{\frac{\pi}{2}}\cdot
x^{\frac{3}{2}}K_{\frac{3}{2}}(x)=\frac{\pi}{2}(x+1)e^{-x},\end{equation} thus from (\ref{eq.34k}) for $\nu=\frac{3}{2}$ we obtain
the following inequality
$$\left|(a+1)e^{-a}-(b+1)e^{-b}\right|\geq |a-b|\sqrt{ab}e^{-\frac{a+b}{2}} \ \ \mbox{for all}\
a,b\geq 0.
$$

\section{\bf Proofs of the main results}
\setcounter{equation}{0}

\begin{proof}[\bf Proof of Theorem \ref{th1}]
(a) \& (b) First observe that for $\nu>-1$ the power series
$$\mathcal{I}'_{\nu}(x)=\sum_{n\geq 1}\frac{\left(\frac{1}{4}\right)^{n-1}}{2{(\nu+1)}_n(n-1)!}\cdot x^{2n-1}$$
has only positive coefficients, and thus it makes sense to study
the log-concavity of $x\mapsto\mathcal{I}_{\nu}'(x).$ Moreover observe that $x\mapsto\mathcal{I}_{\nu}'(x)$ is strictly log-concave on $(0,\infty)$ if and only if $x\mapsto
\mathcal{I}_{\nu}''(x)/\mathcal{I}_{\nu}'(x)$ is strictly decreasing on $(0,\infty).$ By using the recurrence formula \cite[p. 79]{watson} $\left[x^{-\nu}I_{\nu}(x)\right]'=x^{-\nu}I_{\nu+1}(x)$ we have
\begin{equation}\label{eq.23}
\mathcal{I}_{\nu}'(x)=2^{\nu}\Gamma(\nu+1)x^{-\nu}I_{\nu+1}(x).\end{equation}
On the other hand application of \cite[p. 79]{watson} $xI_{\nu+1}(x)=xI_{\nu}'(x)-{\nu}I_{\nu}(x)$ and
(\ref{eq.23}) yields
\begin{equation}\label{eq.27}\mathcal{I}_{\nu}''(x)=2^{\nu}\Gamma(\nu+1)x^{-\nu}\left[I_{\nu+2}(x)+\frac{1}{x}I_{\nu+1}(x)\right].\end{equation}
By using the notation $R_{\nu}(x)=I_{\nu+1}(x)/I_{\nu}(x),$ in view of \eqref{eq.23} and \eqref{eq.27}, we obtain that $x\mapsto
\mathcal{I}_{\nu}''(x)/\mathcal{I}_{\nu}'(x)$ is strictly decreasing on $(0,\infty)$ for all $\nu\in\left(-1,-\frac{1}{2}\right]$ if and only if the function $x\mapsto Q_{\nu+1}(x),$ where $Q_{\nu}(x)={1}/{x}+R_{\nu}(x),$ is strictly decreasing on $(0,\infty)$ for all $\nu\in\left(-1,-\frac{1}{2}\right].$ To prove this we need to show that $Q_{\nu}$ is strictly decreasing on $(0,\infty)$ for all $\nu\in\left(0,\frac{1}{2}\right].$ However, since the function
$$x\mapsto Q_{\frac{1}{2}}(x)=R_{\frac{1}{2}}(x)+{1}/{x}={I_{\frac{3}{2}}(x)}/{I_{\frac{1}{2}}(x)}+{1}/{x}=\coth x$$
is strictly decreasing on $(0,\infty),$ we can suppose that $\nu\in\left(0,\frac{1}{2}\right).$ By using the fact (see \cite[p. 446]{yuan}) that $R_{\nu}(x)\to0$ as $x\to 0,$ we have that $Q_{\nu}(x)\to\infty$ as $x\to0.$ On the other hand, by using the asymptotic expansion of $I_{\nu}(x)$ for large values of $x$ it can be seen (see \cite[p. 446]{yuan}) that $R_{\nu}(x)\to1$ as $x\to\infty,$ and consequently $Q_{\nu}(x)\to1$ as $x\to\infty.$ Thus, for small values of $x$ the function $Q_{\nu}$ is strictly decreasing. Moreover, we can show that $Q_{\nu}$ is strictly decreasing on $(0,\sqrt{2\nu+4}),$ and consequently on $(0,2),$ for all $\nu>0.$ To prove this, first observe that $Q_{\nu}(x)$ can be rewritten as
$$Q_{\nu}(x)=\frac{I_{\nu}'(x)}{I_{\nu}(x)}+\frac{1-\nu}{x}$$ and by using twice the inequality \cite[p. 526]{segura}
$$y_{\nu}(x)=\frac{xI_{\nu}'(x)}{I_{\nu}(x)}>\nu+\frac{x^2}{\nu+\frac{1}{2}+\sqrt{x^2+\left(\nu+\frac{3}{2}\right)^2}},$$
where $\nu>0,$ $x>0,$ we obtain that
\begin{align*}x^2Q_{\nu}'(x)&=xy_{\nu}'(x)-y_{\nu}(x)+\nu-1=x^2+\nu^2-y_{\nu}^2(x)-y_{\nu}(x)+\nu-1\\&<
-\left(\frac{2\nu+1}{\nu+\frac{1}{2}+\sqrt{x^2+\left(\nu+\frac{3}{2}\right)^2}}-1\right)x^2-
\frac{x^4}{\left(\nu+\frac{1}{2}+\sqrt{x^2+\left(\nu+\frac{3}{2}\right)^2}\right)^2}-1\\
&=\frac{2(\nu+1)x^2}{\left(\nu+\frac{1}{2}+\sqrt{x^2+\left(\nu+\frac{3}{2}\right)^2}\right)^2}-1=
\frac{\sqrt{x^2+\left(\nu+\frac{3}{2}\right)^2}-\left(\nu+\frac{1}{2}\right)}{\sqrt{x^2+\left(\nu+\frac{3}{2}\right)^2}+\left(\nu+\frac{1}{2}\right)}
\cdot\frac{2(\nu+1)x^2}{2(\nu+1)+x^2}-1\equiv w_{\nu}(x).
\end{align*}
Now taking into account that for $\nu\geq-\frac{1}{2}$ both functions
$$t\mapsto \frac{\sqrt{t+\left(\nu+\frac{3}{2}\right)^2}-\left(\nu+\frac{1}{2}\right)}{\sqrt{t+\left(\nu+\frac{3}{2}\right)^2}+\left(\nu+\frac{1}{2}\right)}\ \ \ \mbox{and}\ \ \ t\mapsto \frac{2(\nu+1)t}{2(\nu+1)+t}$$
are increasing on $(0,\infty),$ it follows that their product is also increasing and this in turn implies that for $\nu>0$ and $x\in(0,\sqrt{2\nu+4})$ we have
$$w_{\nu}(x)<\frac{\sqrt{2\nu+4+\left(\nu+\frac{3}{2}\right)^2}-\left(\nu+\frac{1}{2}\right)}{\sqrt{2\nu+4+\left(\nu+\frac{3}{2}\right)^2}
+\left(\nu+\frac{1}{2}\right)}\cdot\frac{2(\nu+1)(2\nu+4)}{2(\nu+1)+2\nu+4}-1=\varphi(\nu)<\lim_{\nu\to\infty}\varphi(\nu)=0.$$
Here we used that the function $\varphi$ is increasing on $(0,\infty).$ Summarizing, we have proved that for $\nu>0$ and $x\in(0,\sqrt{2\nu+4})$ the function $Q_{\nu}$ satisfies $x^2Q_{\nu}'(x)<0,$ which implies that $Q_{\nu}$ is strictly decreasing on $(0,2)$ for $\nu>0.$ Thus the first extreme of the function $Q_{\nu},$ if any, is a minimum. Now, since $R_{\nu}(x)$ satisfies the differential equation
$$u'(x)=1-u^2(x)-\frac{2\nu+1}{x}u(x),$$ it follows that $Q_{\nu}(x)$ satisfies the differential equation
$$v'(x)=1-v^2(x)-\frac{2\nu-1}{x}v(x)+\frac{2\nu-1}{x^2}.$$ Differentiating both sides of the above relation we get
$$v''(x)=\frac{2\nu-1}{x^2}\left[v(x)-\frac{2}{x}\right]-\left[2v(x)+\frac{2\nu-1}{x}\right]v'(x).$$
This shows that if $Q_{\nu}'(x_{\star})=0$ for some $x_{\star}\geq2,$ then for that $x_{\star}$ we have
$$Q_{\nu}''(x_{\star})=\frac{2\nu-1}{x_{\star}^2}\left[Q_{\nu}(x_{\star})-\frac{2}{x_{\star}}\right]=\frac{2\nu-1}{x_{\star}^2}\left[R_{\nu}(x_{\star})-\frac{1}{x_{\star}}\right].$$
On the other hand, by using the inequality \cite[p. 241]{amos}
$$R_{\nu}(x)>\frac{x}{\sqrt{x^2+\left(\nu+1\right)^2}+\nu+1},$$
which holds for $\nu\geq 0$ and $x>0,$ we obtain that for $\nu\in\left(0,\frac{1}{2}\right)$ and $x\geq2$ we have
$$R_{\nu}(x)>\frac{x}{\sqrt{x^2+\frac{9}{4}}+\frac{3}{2}}\geq\frac{1}{x}.$$
This implies that $Q_{\nu}''(x_{\star})<0.$ But, this is in contradiction with the fact that the first extreme is a minimum, and consequently the derivative of $Q_{\nu}$ does not vanish when $\nu\in\left(0,\frac{1}{2}\right).$ This in turn implies that $Q_{\nu}$ is strictly decreasing on $(0,\infty)$ for $\nu\in\left(0,\frac{1}{2}\right),$ that is, $Q_{\nu+1}$ is strictly decreasing on $(0,\infty)$ for $\nu\in\left(-1,-\frac{1}{2}\right).$ This is equivalent to the fact that if $\nu\in\left(-1,-\frac{1}{2}\right),$ then the function $x\mapsto\mathcal{I}_{\nu}'(x)$ is strictly log-concave on $(0,\infty).$ On the other hand, since $Q_{\nu}$ is strictly decreasing on $\left(0,\sqrt{2(\nu+2)}\right)$ for $\nu>0,$ we get that $Q_{\nu+1}$ is strictly decreasing on $\left(0,\sqrt{2(\nu+3)}\right)$ for $\nu>-1,$ and this is equivalent to the fact that if $\nu>-1,$ then the function $x\mapsto\mathcal{I}_{\nu}'(x)$ is strictly log-concave on $\left(0,\sqrt{2(\nu+3)}\right).$

(c) By using the above things it is clear that we need to show that for $\nu>\frac{1}{2}$ there exists an $x_{\nu-1}>\sqrt{2(\nu+2)}$ such that $Q_{\nu}$ is strictly decreasing on $(0,x_{\nu-1})$ and strictly increasing on $(x_{\nu-1},\infty).$ For this we note that \cite[p. 276]{gronwall}
$$y_{\nu}(x)-x\sim -\frac{1}{2}+\frac{4\nu^2-1}{8x}-\dots,$$
which holds for large values of $x$ and fixed $\nu.$ From this we obtain that for $\nu>\frac{1}{2}$ and large values of $x$ we have
$y_{\nu}(x)-x<\nu-1,$ or equivalently $Q_{\nu}(x)<1.$ Since $Q_{\nu}$ is strictly decreasing on $\left(0,\sqrt{2\nu+4}\right)$ and $Q_{\nu}(x)\to1$ as $x\to\infty,$ it follows that there must be a point at which $Q_{\nu}'(x)$ vanishes. In other words, there exists an $x_{\nu-1}>\sqrt{2(\nu+2)}$ such that $Q_{\nu}$ is strictly decreasing on $(0,x_{\nu-1})$ and strictly increasing on $(x_{\nu-1},\infty).$ Moreover, $Q_{\nu}'$ changes sign only once, otherwise we will have two local minima between which there must be a local maximum, which is impossible, since in that stationary point $x^{\star}$ we would have
$$Q_{\nu}''(x^{\star})=\frac{2\nu-1}{(x^{\star})^2}\left[R_{\nu}(x^{\star})-\frac{1}{x^{\star}}\right]>0,$$
according to the above discussion.
\end{proof}

\begin{proof}[\bf Proof of Corollary \ref{th2}]
Without loss of generality it is enough to show that
inequality (\ref{eq.34}) holds for $b>a>0.$ Observe that the
function $x\mapsto\mathcal{I}_{\nu}'(x)$ is convex on $[a,b],$ since
as power series has only positive coefficients. Now taking into
account that from Theorem \ref{th1} the function $x\mapsto\mathcal{I}_{\nu}'(x)$ is strictly log-concave on $(0,\infty),$
from the Hermite-Hadamard inequality \cite[p.
14]{mi} applied to the function $x\mapsto\mathcal{I}_{\nu}'(x)$ we
have
$$\frac{\mathcal{I}_{\nu}(b)-\mathcal{I}_{\nu}(a)}{b-a}=\frac{1}{b-a}\int_a^b\mathcal{I}_{\nu}'(x)dx\geq
\mathcal{I}'_{\nu}\left(\frac{a+b}{2}\right)\geq
\sqrt{\mathcal{I}_{\nu}'(a)\mathcal{I}_{\nu}'(b)},$$ which completes the proof.

Another proof of Corollary \ref{th2} is as follows. Let $f:[a,b]\to\mathbb{R}$ be such that $f':[a,b]\to(0,\infty)$ is log-concave. Then the inequality \cite[p. 242]{corput}
\begin{align}\label{fcorput}
\log\left(\frac{1}{b-a}\int_a^bf'(x)dx\right)\geq \frac{1}{b-a}\int_a^b\log f'(x)dx\geq\frac{\log f'(a)+\log f'(b)}{2}
\end{align}
implies that
$$\log \frac{f(b)-f(a)}{b-a}\geq \log\sqrt{f'(a)f'(b)}$$
is valid. Applying \eqref{fcorput} for the function $f=\mathcal{I}_{\nu},$ which has the property that $x\mapsto\mathcal{I}_{\nu}'(x)$ is strictly log-concave on $(0,\infty),$
we arrive at \eqref{eq.34}.
\end{proof}

\begin{proof}[\bf Proof of Theorem \ref{th1bessel}] By using the recurrence relation \cite[p. 45]{watson} $\left[x^{-\nu}J_{\nu}(x)\right]'=-x^{-\nu}J_{\nu+1}(x)$ we get
$$\mathcal{J}_{\nu}'(x)=-2^{\nu}\Gamma(\nu+1)x^{-\nu}J_{\nu+1}(x),$$
$$\mathcal{J}_{\nu}''(x)=2^{\nu}\Gamma(\nu+1)x^{-\nu}\left[J_{\nu+2}(x)-\frac{1}{x}J_{\nu+1}(x)\right],$$
and thus the function $x\mapsto -\mathcal{J}_{\nu}'(x),$ which maps $(0,j_{\nu+1,1})$ into $(0,\infty),$ is strictly log-concave if and only if the function $x\mapsto {1}/{x}-J_{\nu+2}(x)/J_{\nu+1}(x)$ is strictly decreasing on $(0,j_{\nu+1,1}).$ But, using the well-known Mittag-Leffler expansion
$$\frac{J_{\nu+1}(x)}{J_{\nu}(x)}=\sum_{n\geq 1}\frac{2x}{j_{\nu,n}^2-x^2},$$
we obtain that
$$\left[\frac{J_{\nu+1}(x)}{J_{\nu}(x)}\right]'=2\sum_{n\geq 1}\frac{j_{\nu,n}^2+x^2}{(j_{\nu,n}^2-x^2)^2}>0$$
for all $\nu>-1,$ $x\neq j_{\nu,n},$ $n\in\{1,2,\dots\}.$ Note that this can be proved also by using the Neumann's formula for the product of two Bessel functions of the first kind with different orders, see \cite[Lemma 2.5]{im}. In particular, the above result means that the function $x\mapsto J_{\nu+1}(x)/J_{\nu}(x)$ is strictly increasing on $(0,j_{\nu,1})$ for $\nu>-1.$ Thus, for $\nu>-2$ the function $x\mapsto {1}/{x}-J_{\nu+2}(x)/J_{\nu+1}(x)$ is strictly decreasing on $(0,j_{\nu+1,1})$ as the sum of two strictly decreasing functions.
\end{proof}

\begin{proof}[\bf Proof of Corollary \ref{th2bessel}] Observe that as in Corollary \ref{th2} it is enough to consider the case when $b>a.$ Then the proof of the result follows from Theorem \ref{th1bessel} and \eqref{fcorput} by choosing for $f$ the function $-\mathcal{J}_{\nu}.$
\end{proof}

\begin{proof}[\bf Proof of Theorem \ref{th3k}]
By using the recurrence relation \cite[p. 79]{watson} $xK_{\nu}'(x)+\nu K_{\nu}(x)=-xK_{\nu-1}(x)$ we obtain that
$$\mathcal{K}_{\nu}'(x)=-2^{\nu-1}\Gamma(\nu)x^{\nu}K_{\nu-1}(x),$$
$$\mathcal{K}_{\nu}''(x)=-2^{\nu-1}\Gamma(\nu)x^{\nu}\left[\frac{\nu}{x}K_{\nu-1}(x)+K_{\nu-1}'(x)\right],$$
and thus $x\mapsto -\mathcal{K}_{\nu}'(x)$ is strictly log-concave on $(0,\infty)$ if and only if the function
$$x\mapsto \frac{\mathcal{K}_{\nu}''(x)}{\mathcal{K}_{\nu}'(x)}=\frac{\nu}{x}+\frac{K_{\nu-1}'(x)}{K_{\nu-1}(x)}=\frac{1}{x}-\frac{K_{\nu-2}(x)}{K_{\nu-1}(x)}$$
is strictly decreasing on $(0,\infty).$ By using \cite[Lemma 2.4]{im} the function $x\mapsto K_{\nu+1}(x)/K_{\nu}(x)$ is strictly decreasing on $(0,\infty)$ for $\nu>-\frac{1}{2},$ which implies that the function $x\mapsto {K_{\nu-2}(x)}/{K_{\nu-1}(x)}$ is strictly increasing on $(0,\infty)$ for $\nu>\frac{3}{2}.$ Consequently, the function $x\mapsto -\mathcal{K}_{\nu}'(x)$ is strictly log-concave on $(0,\infty)$ for $\nu>\frac{3}{2}.$ Finally, observe that by using \eqref{k3/2} we have that $x\mapsto -\mathcal{K}_{\frac{3}{2}}'(x)=\frac{\pi}{2}xe^{-x}$ is also strictly log-concave on $(0,\infty).$ This completes the proof.
\end{proof}

\begin{proof}[\bf Proof of Corollary \ref{th4k}]
Clearly, as in Corollary \ref{th2} it is enough to consider the case when $b>a.$ Then the proof of the result follows from Theorem \ref{th3k} and \eqref{fcorput} by choosing for $f$ the function $-\mathcal{K}_{\nu}.$
\end{proof}

\section{\bf Concluding remarks}
\setcounter{equation}{0}

In this section we would like to comment the main results of this paper.

\begin{enumerate}

\item[\bf 1.] We note that in \cite[Theorem 3.6.3]{bariczPhD} (see also \cite[Theorem 3.25]{barspringer}) it was proved that the function $x \mapsto \mathcal{I}_{\nu}'(x)$ is strictly log-concave on $\left(0,\sqrt{2(\nu+2)}\right)$ for $\nu> -1.$ Moreover, it was shown in \cite[Theorem 3.6.4]{bariczPhD} (see also \cite[Theorem 3.26]{barspringer}) that the inequality \eqref{eq.34} is valid for all $\nu>-1$ and $a,b\in\left(0,\sqrt{2(\nu+2)}\right).$ Theorem \ref{th1} and Corollary \ref{th2} complement and improve the above mentioned results.

\item[\bf 2.] We also note that (see \cite[Theorem 3.6.3]{bariczPhD} or \cite[Theorem 3.25]{barspringer}) for $\nu>-1$ the function $x \mapsto \mathcal{I}_{\nu}(x)$ is strictly log-convex on $\left(0,\sqrt{2(\nu+1)(\nu+2)}\right).$ Since the function $x \mapsto \log\mathcal{I}_{\nu}(x)$ is even it follows that in fact $x \mapsto \mathcal{I}_{\nu}(x)$ is strictly log-convex on $\left(-\sqrt{2(\nu+1)(\nu+2)},\sqrt{2(\nu+1)(\nu+2)}\right)$ for $\nu>-1.$ Moreover, it was shown in \cite[Theorem 4]{barexpo} that the function $x \mapsto \mathcal{I}_{\nu}(x)$ is strictly log-convex on $[-j_{\nu,1},j_{\nu,1}]$ for $\nu>-1,$ where $j_{\nu,1}$ is the first positive zero of the Bessel function of the first kind $J_{\nu}.$ Next, we show that the function $x \mapsto \mathcal{I}_{\nu}(x)$ cannot be strictly log-convex on the whole $\mathbb{R}$ for $\nu>-1.$ When $\nu\geq -\frac{1}{2}$ the strict log-convexity of $x \mapsto \mathcal{I}_{\nu}(x)$ is certainly valid, but in the case when $\nu\in\left(-1,-\frac{1}{2}\right)$ there exists a $z_{\nu}\geq j_{\nu,1}$ such that the function $x \mapsto \mathcal{I}_{\nu}(x)$ is strictly log-convex on $(0,z_{\nu})$ and strictly log-concave on $(x_{\nu},\infty).$ For this observe that the function $x \mapsto \mathcal{I}_{\nu}(x)$ is strictly log-convex (log-concave) if and only if $x\mapsto R_{\nu}(x)$ is strictly increasing (decreasing). But, it was shown in \cite[p. 446]{yuan} that when $\nu\in\left(-1,-\frac{1}{2}\right)$ the function $x\mapsto R_{\nu}(x)$ on $(0,\infty)$ is increasing first to reach a maximum and then decreasing. This means that there exists a $z_{\nu}$ such that $x\mapsto R_{\nu}(x)$ is strictly increasing on $(0,z_{\nu})$ and strictly decreasing on $(z_{\nu},\infty),$ as we required. Now, since $x \mapsto \mathcal{I}_{\nu}(x)$ is strictly log-convex on $(0,j_{\nu,1}]$ we clearly have that $z_{\nu}\geq j_{\nu,1}.$

\item[\bf 3.] We would like to mention that the idea of the first proof of \eqref{eq.34} does not work in the case of \eqref{eq.34bessel} because in view of
$\mathcal{J}_{-\frac{1}{2}}(x)=\cos x,$ the function
$x\mapsto -\mathcal{J}'_{-\frac{1}{2}}(x)=\sin x$ is concave and log-concave and consequently
$$\mathcal{J}_{-\frac{1}{2}}'\left(\frac{a+b}{2}\right)\leq \frac{1}{b-a}\int_a^b\mathcal{J}_{-\frac{1}{2}}'(x)dx=\frac{\mathcal{J}_{-\frac{1}{2}}(b)-\mathcal{J}_{-\frac{1}{2}}(a)}{b-a},\ \ \ -\mathcal{J}_{-\frac{1}{2}}'\left(\frac{a+b}{2}\right)\geq \sqrt{\mathcal{J}_{-\frac{1}{2}}'(a)\mathcal{J}_{-\frac{1}{2}}'(b)}$$
and these together does not imply the van der Corput's inequality \eqref{eq.32sine}. However, as we have seen in the proof of Corollary \ref{th2bessel} the inequality \eqref{fcorput} applied to the function $-\mathcal{J}_{\nu}$ implies the van der Corput inequality \eqref{eq.34bessel}. In other words, there no need to assume for $f'$ to be convex, its log-concavity is quite enough, according to \eqref{fcorput}.

\item[\bf 4.] Finally, we would like to mention that functions which have log-concave derivatives occur for example in probability theory. Namely, many probability density functions are log-concave (see for example \cite{bargeo}), and thus applying the inequality \eqref{fcorput} for the cumulative distribution functions we can get van der Corput inequalities. For example, consider the probability density and cumulative distribution function of the standard normal distribution
$$\varphi(x)=\frac{1}{\sqrt{2\pi}}e^{-\frac{x^2}{2}}\ \ \ \mbox{and}\ \ \ \Phi(x)=\frac{1}{\sqrt{2\pi}}\int_{-\infty}^xe^{-\frac{t^2}{2}}dt.$$
The function $x\mapsto \Phi'(x)=\varphi(x)$ is log-concave on $\mathbb{R},$ and applying the inequality \eqref{fcorput} for the function $f=\Phi$ we get the following van der Corput inequality for all real $a$ and $b$
$$|\Phi(a)-\Phi(b)|\geq |a-b|\sqrt{\varphi(a)\varphi(b)}.$$
Another example is the gamma distribution which has support $(0,\infty),$ probability density and cumulative distribution function as follows
$$f(x)=\frac{x^{\alpha-1}e^{-x}}{\Gamma(\alpha)}\ \ \ \mbox{and} \ \ \  F(x)=\frac{\gamma(\alpha,x)}{\Gamma(\alpha)}=\frac{1}{\Gamma(\alpha)}\int_0^x t^{\alpha-1}e^{-t}dt,$$
where $\gamma(\alpha,\cdot)$ is the lower incomplete gamma function and $\alpha>0$ is the shape parameter. Since the probability density function of the gamma distribution is log-concave on $(0,\infty)$ for all $\alpha\geq 1$ (see for example \cite[p. 192]{bargeo}), applying \eqref{fcorput} for the function $F$ we get the following van der Corput inequality
$$|\gamma(\alpha,a)-\gamma(\alpha,b)|\geq |a-b|(ab)^{\frac{\alpha-1}{2}}e^{-\frac{a+b}{2}},$$
where $a,b\geq 0$ and $\alpha\geq 1.$
\end{enumerate}

\subsection*{Acknowledgement} The work of \'A. Baricz was supported by the J\'anos Bolyai Research Scholarship of
the Hungarian Academy of Sciences. This author is very grateful to Aingeru Fern\'andez from University of the Basque Country, Spain, for pointing out an error in the proof of Theorem 1 in an earlier version of the paper.


\begin{thebibliography}{}

\bibitem{amos} D.E. Amos, Computation of modified Bessel functions and their ratios, Math. Comp. 28 (1974) 239--251.

\bibitem{barexpo} \'A. Baricz, Functional inequalities involving Bessel and modified
Bessel functions of the first kind, Expo. Math. 26 (2008) 279--293.

\bibitem{bariczPhD} \'A. Baricz, Generalized Bessel Functions of the First Kind, PhD Thesis, Babe\c{s}-Bolyai
University, Cluj-Napoca, 2008.

\bibitem{barspringer} \'A. Baricz, Generalized Bessel functions of the First Kind, Lecture Notes in Mathematics, vol. 1994, Springer, Berlin, 2010.

\bibitem{bargeo} \'A. Baricz, Geometrically concave univariate distributions, J. Math. Anal. Appl. 363 (2010) 182--196.

\bibitem{bane1} \'A. Baricz, E. Neuman, Inequalities involving modified Bessel
functions of the first kind II, J. Math. Anal. Appl. 332 (2007) 265--271.

\bibitem{corput} J.G. van der Corput, Problem 95, Wisk. Opgaven 16 (1935) 242--243.

\bibitem{GLP}
C. Giordano, A. Laforgia, J. Pe\v{c}ari\'c, Supplement to
known inequalities for some special functions, J. Math. Anal.
Appl. 200 (1996) 34--41.

\bibitem{gronwall}
T.H. Gronwall, An inequality for the Bessel functions of the first kind with imaginary argument, Ann. of Math. 33(2)
(1932) 275--278.

\bibitem{MI}
M.E.H. Ismail, Remarks on a paper by Giordano, Laforgia,
and Pec\v{a}ri\'c, J. Math. Anal. Appl. 211 (1997) 621--625.

\bibitem{im}
M.E.H. Ismail, M.E. Muldoon, Monotonicity of the zeros of a cross-product of
Bessel functions, SIAM J. Math. Anal. 9(4) (1978) 759--767.

\bibitem{mi} {D.S. Mitrinovi\'c,} Analytic
Inequalities, Springer-Verlag, Berlin, 1970.

\bibitem{neuman} E. Neuman, Inequalities involving modified Bessel
functions of the first kind, J. Math. Anal. Appl. 171 (1992) 532--536.

\bibitem{segura} J. Segura, Bounds for ratios of modified Bessel functions and associated Tur\'an-type
inequalities, J. Math. Anal. Appl. 374 (2011) 516--528.

\bibitem{yuan} L. Yuan, J.D. Kalbfleisch, On the Bessel distribution and related problems, Ann. Inst. Statist. Math. 52 (2000) 438--447.

\bibitem{watson} G.N. Watson, A Treatise on the Theory of Bessel
Functions, Cambridge University Press, Cambridge, 1944.

\end{thebibliography}
\end{document}